\documentclass[preprint]{elsarticle}

\biboptions{sort&compress}

\usepackage{amsmath,amssymb}
\usepackage{graphicx}

\usepackage[utf8]{inputenc}
\usepackage{mathtools}  
\usepackage{amsthm}
\usepackage[colorlinks]{hyperref}
\usepackage[arrowdel]{physics}

\usepackage{cleveref}

\usepackage[dvipsnames]{xcolor}

\usepackage{bm}
\let\vec\mathbf

\newcommand{\I}{\mathcal{I}}

\newcommand{\vx}{\vec{x}}

\newcommand{\vecu}{\vec{u}}
\newcommand{\vmu}{\boldsymbol{\mu}}
\newcommand{\vbeta}{\boldsymbol{\beta}}

\newcommand{\vecl}{\boldsymbol{\lambda}}
\newcommand{\J}{\mathcal{J}}

\newcommand{\mL}{\mathcal{L}}
\newcommand{\mH}{\mathcal{H}}

\newcommand{\mub}{\bar{\mu}}

\newcommand{\de}{\partial}
\newcommand{\Gl}{G_{\vecl}}
\newcommand{\dq}{\delta q}
\newcommand{\del}{\de_{\vecl}}
\newcommand{\deu}{\de_{\vecu}}
\newcommand{\deq}{\de_{q}}

\newcommand{\R}{\mathbb{R}}
\newcommand{\model}{\mathcal{M}}

\newcommand{\nhat}{\hat{\vec{n}}}

\newcommand{\citeg}[1]{\cite[e.g.][]{#1}}

\begin{document}

\begin{frontmatter}

\title{Exact vs approximate second-order derivatives in vertically-integrated ice sheet models}

\author[inst1]{Trystan Surawy-Stepney\corref{c1}}
\author[inst2]{Stephen L. Cornford}
\address[inst1]{University of Leeds, Leeds, United Kingdom}
\cortext[c1]{Corresponding author: t.surawystepney@leeds.ac.uk}
\address[inst2]{University of Bristol, Bristol, United Kingdom}

\begin{abstract}
Second order derivatives of model outputs with respect to input parameters are key to several applications in ice sheet modelling. For example, the ability to compute Hessian-vector products broadens the list of available optimisation methods, and facilitates certain kinds of parametric uncertainty quantification. Some modern ice sheet models are built on frameworks supporting algorithmic differentiation (AD), allowing for the computation of higher order derivatives with relative ease. However, many of our most widely-used models are not. A natural alternative might be to follow common practise in first order gradient computation and construct an approximate second-order adjoint model at the PDE level, which neglects the nonlinear dependence of ice viscosity on velocity. Here, we present such a model for the shallow-stream approximation allowing one to compute approximate second-order derivatives, and compare with full second-order derivates found using AD. We find that this produces Hessian-vector products that are superficially similar to those computed via AD. However, an analysis of the spectral decomposition of the Hessians calculated in each way reveals that the subspaces spanned by their eigenvectors diverge after the leading 4 modes, though divergence does not accelerate after this. We conclude that the utility of the approximate Hessian is case-dependent, and a full Hessian, likely computed using AD, should be used where high fidelity is required above very low rank.
\end{abstract}

\begin{keyword}
algorithmic differentiation \sep adjoints \sep uncertainty quantification \sep inverse problems \sep ice sheets \sep Hessian \sep differentiable programming
\end{keyword}

\end{frontmatter}

\section{Introduction}

The Earth's ice sheets, both present and past, are treated by the majority of ice sheet models as continuum materials that undergo viscous deformation described by the Stokes equations. Numerical models typically solve this set of momentum-balance equations, or some lower-dimensional approximation, for ice flow velocity in order to determine how ice is redistributed over time. In addition, many contemporary models contain machinery to compute gradients of model outputs with respect to certain input parameters - driven largely by a need to infer unknown basal or rheological properties of the ice using satellite-derived observations, for example, for the purpose of model initialisation \cite{seroussi_initmip}. For most models, this kind of inference is done considering the elliptic momentum balance problem in isolation (the case we focus on here) \citeg{MacAyeal_1993, Rommelaere_MacAyeal_1997, vieli2003application, joughin2004RossIPs, Larour2005Inverse, Joughin_PIGIPs2009, morlighem2010, cornford2015century}, though some can accommodate time dependency as well \cite{goldbergheimbach2013, Larour2014, Koziol2021}.\\

There are two main methods for computing these gradients. The first is to rely on algorithmic differentiation (AD), should the model be built on a framework that allows it. This technology, used widely in the modelling of Earth system dynamics since the 1990s \citeg{WhatIsanAdjointModel, marshallFVNS1997, Giering1998, heimback2002, Heimbach_Bugnion_2009}, in theory, facilitates the computation of derivatives of complex functions expressed in code with little developer overhead. Programmatic operations are accompanied by rules defining the action of derivatives of those operations on tangents (known as \emph{forward mode} AD) or cotangents (known as \emph{reverse} or \emph{adjoint mode} AD) to the parameter space. These rules can be derived before or at runtime, and there are a wide variety of techniques for doing so \cite{griewank2008evaluating}. Various widely-used ice sheet models include some variant of this technology \cite{goldbergheimbach2013, Hück02112018, icepack2021}.\\

The alternative to AD, in the computation of first-order derivatives, is to use PDE-level adjoint methods. Here, gradients are found by solving a linear PDE for a set of adjoint variables which keep track of perturbations to the model outputs that maintain the validity of the model equations under input perturbations \cite{lions1971optimal}. To remove some of the algebraic complexity, the first-order adjoint methods implemented in many models assume a linear rheology \citeg{MacAyeal_1993, Rommelaere_MacAyeal_1997,morlighem2010,cornford2015century}, resulting in a self-adjoint set of model equations and approximate gradients. Of course, unless the approximate gradients are (roughly) parallel to the actual gradients, optimisation algorithms that depend on them will likely converge to different solutions. The use of the self-adjoint approximation has been found to produce inaccurate or low resolution solutions for scalar basal slipperiness parameters (widely used in ice sheet modelling) in various experiments involving synthetic data \cite{goldberg_sergienko_2011, martin_monier_2014}. However, the difference in performance when applied to real settings is thought to be much smaller \cite{goldberg_sergienko_2011}, and has been shown, for some pan-Antarctic inverse problems, to be virtually zero \cite{morlighem_ant_2013}. This latter viewpoint has come to dominate over the last decade, and the widespread use of the linear-viscosity approximation in constructing first-order adjoints reflects the rarity of major errors arising from it.\\

Beyond the computation of gradients, there is growing interest in second-order derivative information among ice sheet modellers. For example, Hessian-vector products enable certain powerful optimisation methods (e.g. Newton's method), and are essential for many Bayesian frameworks for parametric uncertainty quantification \cite{thacker1988, thacker_roleOfHessian_1989, MacAyeal_1993, petra2014Bayes, bea_2023_paper, bea_ice_stream_sensitivity}. To date, codes involving second-order derivatives in the ice sheet modelling community have been based on finite-element libraries (like FEnICs/Firedrake) employing a particular type of AD applicable to that domain \cite{icepack2021, Koziol2021}. \\

Models lacking AD can, once again, use PDE-level adjoint methods to compute these second-order derivatives - at higher accuracy and efficiency than methods based, at least in part, on finite differences (e.g. \cite{MacAyeal_1993}). Extending the logic of the first-order method described above, second-order derivatives can be computed by constructing two further adjoint variables to keep track of: 1) second-order changes to the state variable and 2) first-order changes to the first-order adjoint variable (see, e.g. \cite{CACUCI2015687, Cacuci01092016}). This incurs two additional linear solves. In principle, this can be repeated ad infinitum, computing $n(n+1)/2$ adjoint variables, using $n(n+1)/2$ linear solves, in pursuit of the $n^{\rm th}$-order derivative. The barrier to implementation is not so much the computational overhead, but the algebraic one in deriving the equations each adjoint variable must obey. However, pockets of the ice sheet modelling community have utilised complete second-order adjoint methods, formulated using variational methods, to good effect: in performing inverse problems with finite element full-Stokes solvers, demonstrating the efficacy of the method in improving the efficiency of inverse solvers and its utility in uncertainty quantification \cite{Petra_Zhu_2012, ISAAC2015348, zhu_geothermal_2016}.\\

Here, we derive a PDE-level second-order adjoint formulation for the shallow-stream approximation (SSA) to the momentum balance, using a linear viscosity approximation to mirror common practice in first-order adjoint methods. We label this the second-order self-adjoint (SOSA) method. The resulting formulation can, therefore, be straightforwardly integrated into a variety of existing SSA (or SSA-like) models and used to compute approximate Hessian-vector products. We implement this model in a finite volume code written in the numerical computation library JAX \cite{jax2018github}. The model is similar to a uniform-mesh version of the BISICLES ice sheet model \cite{CORNFORD2013529}, but fully differentiable, allowing us to compare the SOSA Hessian with an exact Hessian obtained via AD. We qualitatively compare Hessian-vector products (HVPs) in different cases, and perform spectral analysis on each Hessian. We find that, while HVPs from the two methods appear qualitatively similar, the invariant subspaces spanned by their leading eigenvectors begin to diverge almost immediately. However, this divergence ceases after around 50 modes, suggesting a diminishing importance in the neglected nonlinearity. We judge the full Hessian to be preferable where possible, though there may be situations in which SOSA can provide an adequate approximation, with fidelity that depends on the rank of the approximation.\\

\section{The adjoint formulation}
First, we briefly outline the first- and second-order adjoint models for a generic numerical model, before specifying the equations for the SSA formulation of the momentum balance under the self-adjoint approximation. For a full derivation, see \ref{appsec:foa} and \ref{appsec:soa}.\\

Let $U$ and $Q$ be Hilbert spaces with the standard L2 inner product $\langle\cdot,\cdot\rangle$. Let us assume a model $\model$ with input parameters $q\in Q$ and outputs $\vecu\in U$, such that $\model: q\mapsto \vecu$, governed by the system of partial differential equations $G(\vecu, q)=0$, and a functional of interest $\J:U\times Q\to \R$. This could be, for example, the total grounding line flux over some portion of an ice sheet or a cost function comparing modelled and observed ice velocity. The quantities of interest are the gradient $d_q\J:T_Q\to \R$ and the Hessian $d^2_q\J = \mH_q:T_Q\times T_Q\to \R$, where $T_Q$ is the tangent space to $Q$ at the current parameter value.\\

We write the action of the gradient on a perturbation $\delta q$ as $d_q\J[\delta q]=\langle d_q\J, \delta q\rangle$ and the action of the Hessian as $d^2_q\J[\delta q_1, \delta q_2] = \langle d^2\J[\delta q_1], \delta q_2\rangle$. The Hessian-vector-product $d^2_q\J[\delta q_1]$ is itself a gradient operator. The notation here relies on the fact that $T_q\cong Q$, $T_\vecu\cong U$ and the Reisz representation theorem to associate gradients with elements of $U$ and $Q$.\\

\subsection{General first- and second-order adjoint models}

The total variation of the functional $\J$ under a perturbation $\delta q$ has, according to the chain rule, two components:
\begin{equation}\label{dJ}
    d_q\J[\delta q] = (\partial_q\J + \partial_\vecu\J~d_q\vecu)[\delta q].
\end{equation}
So we write $d_q\J = \partial_q\J + \partial_\vecu\J~d_q\vecu$. The adjoint method is a handy trick for calculating this without explicitly computing the Jacobian $d_q\vecu$, which would be prohibitive. Under a small perturbation $\delta q$ that preserves the properties of the original parameter field, the solver will converge just as well, but to a slightly different $\vecu$. Hence:
\begin{align}
d_q\J [\delta q] &= d_q\J[\delta q] + d_q G[\delta q][\vecl]~~~~~~ \forall~\vecl\in U
\end{align}
because $d_q G[\delta q][\vecl] = 0~~~ \forall~\vecl\in U$ .\\

The trick is to choose the adjoint variable $\vecl$ such that this formula contains no reference to the troublesome Jacobian $d_q\vecu$. Following it through (\ref{appsec:foa}), one finds that if $\vecl$ solves the following system of PDEs:
\begin{equation}\label{foa_system}
\partial_\vecu\J + (\partial_\vecu G)^\ast[\vecl] = 0,
\end{equation}
then:
\begin{equation}\label{foa_grad}
d_q\J = \partial_q\J + (\partial_q G)^\ast [\vecl].
\end{equation}
The equations \eqref{foa_system} define the First Order Adjoint (FOA) system; so-called because $(\deu G)^\ast$ is the adjoint of the tangent linear operator. The parameter $\vecl$ encodes the changes to $\vecu$ required to satisfy the model equations under perturbations of the input parameters.\\

Many applications require only that we can compute Hessian-vector products of the form $d^2_q\J [\dq]$, rather than an explicit representation of the full Hessian. The second order adjoint method allows us to do this efficiently by applying the same trick as the first order method, this time defining an additional two adjoint variables $\vmu,\vbeta \in U$. These, respectively, keep track of second order changes to $\vecu$ under perturbations required to maintain the validity of the model equations, and first order changes to the adjoint variable $\vecl$ required to maintain the validity of the FOA equations. The derivation is presented in \ref{appsec:soa}.\\

Let us define:
\begin{equation}
    \Gl = \langle\vecl, G\rangle
\end{equation}
Then the Hessian-vector product can be written:
\begin{equation}\label{generic_hvp}
    d^2_q\J[\dq] = (\deq^2\J + \deq^2\Gl)[\dq] + (\deq G)^\ast[\vbeta] + (\deu\deq\J + (\deq\Gl)^\ast)[\vmu]
\end{equation}
where $\vmu$ and $\vbeta$ solve a set of second-order adjoint (SOA) equations of the form:
\begin{align}
    &\deq\del\Gl[\dq] + \deu\del\Gl[\vmu] = 0\label{mu_eq}\\
    &(\deq\deu\J + \deq\deu\Gl)[\dq] + (\deu G)^\ast[\vbeta] + (\deu^2\J + \deq^2\Gl)[\vmu] = 0.\label{beta_eq}
\end{align}

\subsection{Application to the shallow stream equations}

The FOA and SOA systems given above are generic to any model solving some system of equations expressible as $G(\vecu, q)=0$. In this study, we look specifically at those which solve the shallow stream approximation (SSA) formulation of the momentum balance governing an ice stream \cite{macayeal_1989}. We shall assume without loss of generality that the stress at the base of the ice is a linear function of basal ice velocity. We focus on gradients with respect to the stiffness of the ice, and this is trivially extended to include the slipperiness.\\

Defining the scalar field $\varphi(q)=\phi_0 e^{q}$ and the resistive stress tensor field $H(\vecu) = h\mub\left[\grad\vecu+(\grad\vecu)^\top+2(\grad\cdot\vecu)\I\right]$, the SSA equations read:

\begin{equation}\label{eq:ssa}
    G(q, \vecu) = \grad\cdot[\varphi(q)H(\vecu)] - C\vecu - \rho_i gh\grad s = 0,
\end{equation}

where $h$ is the thickness of the ice, $\vecu$ is the two-dimensional ice velocity, $\mub$ is the vertically-averaged effective viscosity, $C$ is a basal sliding coefficient, $\rho_i$ is the density of ice, $g$ is the acceleration due to gravity, and $s$ is the height of the surface of the ice above sea level.\\

For any $\vecl, G(\vecu, q) \in U$ and $\delta q \in Q$, we can write:
\begin{equation}\label{obvs_identity}
    (\partial_{\vecu,q} G)^\ast[\vecl] = (\partial_{\vecu,q} \langle \vecl, G \rangle).
\end{equation}
Given this, we can write the FOA system as:
\begin{equation}\label{foa_system_ssa}
\partial_\vecu\J + \partial_\vecu \Gl = 0,
\end{equation}
and:
\begin{equation}
d_q\J = \partial_q\J + \partial_q \Gl
\end{equation}
where $\Gl=\langle \vecl,G(\vecu,q) \rangle$, as above.\\

\subsubsection{First order adjoint model}
We consider cases in which the scalar and vector fields of interest and their spatial gradients vanish on the boundaries of the domain. This restricts $U$ and $Q$ from being generic Hilbert spaces. With this, there is a natural equivalence between taking the adjoint of an operator and performing integration by parts. It is easy to ensure this with the finite volume method. When using finite elements, the boundary conditions imposed on $\vecu$ the edges of the domain must also be imposed on the adjoint variables.\\

With this, the first order adjoint system for the SSA equations is:
\begin{equation}\label{foa_ssa}
    \nabla\cdot(\varphi H(\vecl)) - C\vecl + \deu\J = 0.
\end{equation}

With $\vecl$ as the solution to this boundary value problem, the gradient can be written:
\begin{equation}\label{gradient_ssa}
    d_q\J = \deq\J - \nabla\vecl:\varphi H(\vecu).
\end{equation}

\subsubsection{Second order adjoint model}

There are 12 terms to consider in the SOA model and the Hessian-vector product defined in equations (\ref{generic_hvp}-\ref{beta_eq}). However, once more making the assumption of linearity in the viscosity, their form becomes familiar again. The SOA system for the SSA equations is the following set of elliptic PDEs:

\begin{align}
    &\grad\cdot(\varphi\dq H(\vecu)) + \nabla\cdot(\varphi H(\vmu)) - C\vmu = 0\label{mu_eqn}\\
    &\grad\cdot(\varphi\dq H(\vecl)) + \nabla\cdot(\varphi H(\vbeta)) - C\vbeta + (\deu^2\J, \vmu)= 0\label{beta_eqn}.
\end{align}
We have made the assumption here that the functional of interest is separable in functions involving $\vecu$ and $q$ so that $\deq\deu\J=0$.\\

With $\vmu$ and $\vbeta$ found by solving the above, the Hessian vector product is:

\begin{align}\label{hvp_final}
    (\mH,\dq) = (\de^2_q\J,\delta q) - \nabla\vecl:(\varphi\delta q H(\vecu)) -\nabla\vecl:\varphi H(\vmu) -\nabla\vbeta:\varphi H(\vecu).
\end{align}

This shows that, under the self-adjoint approximation, the first- and second-order systems should be solvable by any code built for the SSA equations. We see the same terms repeatedly appearing: the gradient of something that looks like the stress tensor, the linear sliding term and a right-hand side that depends on the form of our functional.

\section{The algorithmic differentiation formulation}

Starting with equation \eqref{dJ}, algorithmic differentiation provides an alternative set of methods for dealing with the term $\deu\J d_q\vecu$. A simple method, implemented, for example, by neural networks, is to unroll the set of elementary operations $\{f_i\}$ that make up the mapping between $q$ and $\vecu$. Imagine $\vecu = f_n\circ f_{n-1}\circ\cdots\circ f_0(q)$, then:
\begin{equation}
    (\deu\J~d_q\vecu)|_q = (\deu\J|_\vecu)(\de_{f_{n}}\vecu|_{f_{n}})(\de_{f_{n-1}}f_{n}|_{f_{n-1}})\cdots(\deq f_0|_{q}),
\end{equation}
(though, in all likelihood, there will be a more complicated computational graph than this structure shows, probably consisting of many branches). In adjoint mode AD (also known as `back-propagation') the series of operations is performed left-to-right. This is significantly more computationally efficient than performing the operations right-to-left when $f_i$ are high-dimensional and the final operation is a functional. However, in a numerical model, the number of operations is huge, and stacking them in this way is not possible.\\

Instead of this, we find a more computationally tractable approach is to implement equations \eqref{foa_system} and \eqref{foa_grad} but use AD to define the action of $(\deu G)^\ast$ and $(\deq G)^\ast$ on adjoint variables $\vecl$ fully - including the nonlinear dependence of ice viscosity on the speed $\vecu$. The function for computing these gradients can then, itself, be differentiated using AD. We use adjoint mode AD for this task, meaning that functions must be specified for the propagation of adjoint variables through a linear solver, but this is again a simple application of a first-order adjoint method.

\section{Implementation}
We make use of a numerical computation library for Python called JAX \cite{jax2018github} to implement the SSA model (eq. \ref{eq:ssa}) and its various adjoints (equations \ref{foa_ssa}-\ref{hvp_final}). With this, we solve the problems on a uniform mesh using the finite volume method. For the forward problem, we use Newton's method, and utilise JAX's AD functionality to compute the Jacobian of the residual equations during each iteration. We use PETSc's python interface \cite{petsc-web-page} to solve linear algebra problems directly using MUMPS \cite{mumps}. To compute the spectral decompositions of the Hessians, we used the Lanczos method using SciPy's sparse.linalg.eigsh function \cite{2020SciPy-NMeth}. Throughout, we use a Glen's flow law exponent of $n=3$.\\

\subsection{Test domains}

\begin{figure}[h!t]
\centering
\includegraphics[width=\linewidth,keepaspectratio]{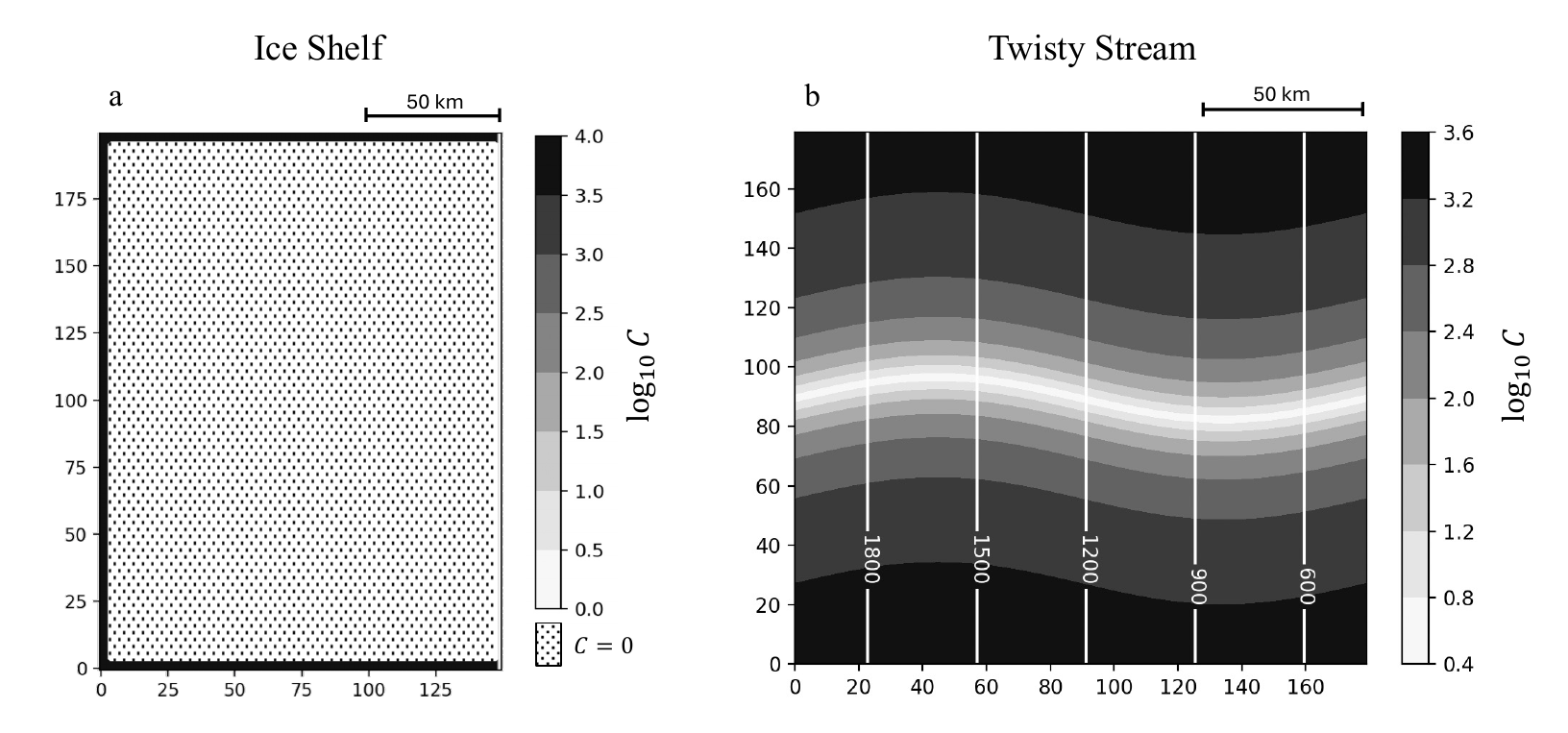}
\caption{Test domains for our experiments. a) The Ice Shelf experiment. Ice is uniform thickness of $500~{\rm m}$ until 2 cells in from the righ-hand-side of the domain, after which thickness is set to $0~{\rm m}$. High values for the slip parameter $C$ are shown in black - represeting grounded ice. The dotted area shows where $C=0$, i.e. floating ice. In the strip of ice-free cells $C$ is set to 1 so that the problem is not underdetermined. b) The Twisty Stream experiment. Ice is set to uniform thickness of $1~{\rm km}$ across the whole domain. The basal slip parameter is shown by the filled contours. White contour lines indicate the surface elevation above sea-level - this decreases uniformly from left to right.}\label{fig:domains}
\end{figure}

To compare the Hessians produced using the AD and SOSA methods, we constructed two synthetic domains: one representing an ice shelf and the other an ice stream (figure \ref{fig:domains}). The first (labelled ``Ice Shelf'' - figure \ref{fig:domains}a) consists of a slab of ice of uniform $500~{\rm m}$ thickness, and a slipperiness coefficient $C$ defined to be zero under most of the ice, and high ($C=10^4$) in a narrow band around three of the sides. The calving front is situated two cells in from the right of the domain. We applied reflection boundary conditions to the velocity at each boundary, though the high $C$-value at three of the boundaries, and the lack of ice at the other, means that the velocity and its gradients are zero there.\\

The other domain (labelled ``Twisty Stream'' - figure \ref{fig:domains}b) is borrowed from \cite{CORNFORD2013529}. It represents an infinitely long, uniform-thickness ($1~{\rm km}$) ice stream following the line of a snaking valley. The slipperiness coefficient is defined as
\[
C = 10^3\times (1 + \epsilon + \sin{\left(\frac{\pi}{2} + 2\pi\frac{y}{R} + m\sin{(2\pi \frac{x}{R})}\right)}
\]
where $\epsilon=5\times10^{-3}$, $R=180~{\rm km}$, and $m=1/4$. The bed is sloping to the right with uniform slope of $0.5^\circ$. To this, we applied periodic boundary conditions on the the left- and right-hand boundaries, and reflection boundary conditions at the top and bottom.\\

The speed solutions for these domains are shown in figure \ref{fig1}a and \ref{fig1}e, respectively.

\section{Results}
\subsection{Direct comparisons}

\begin{figure}[h!t]
\centering
\includegraphics[width=\linewidth,keepaspectratio]{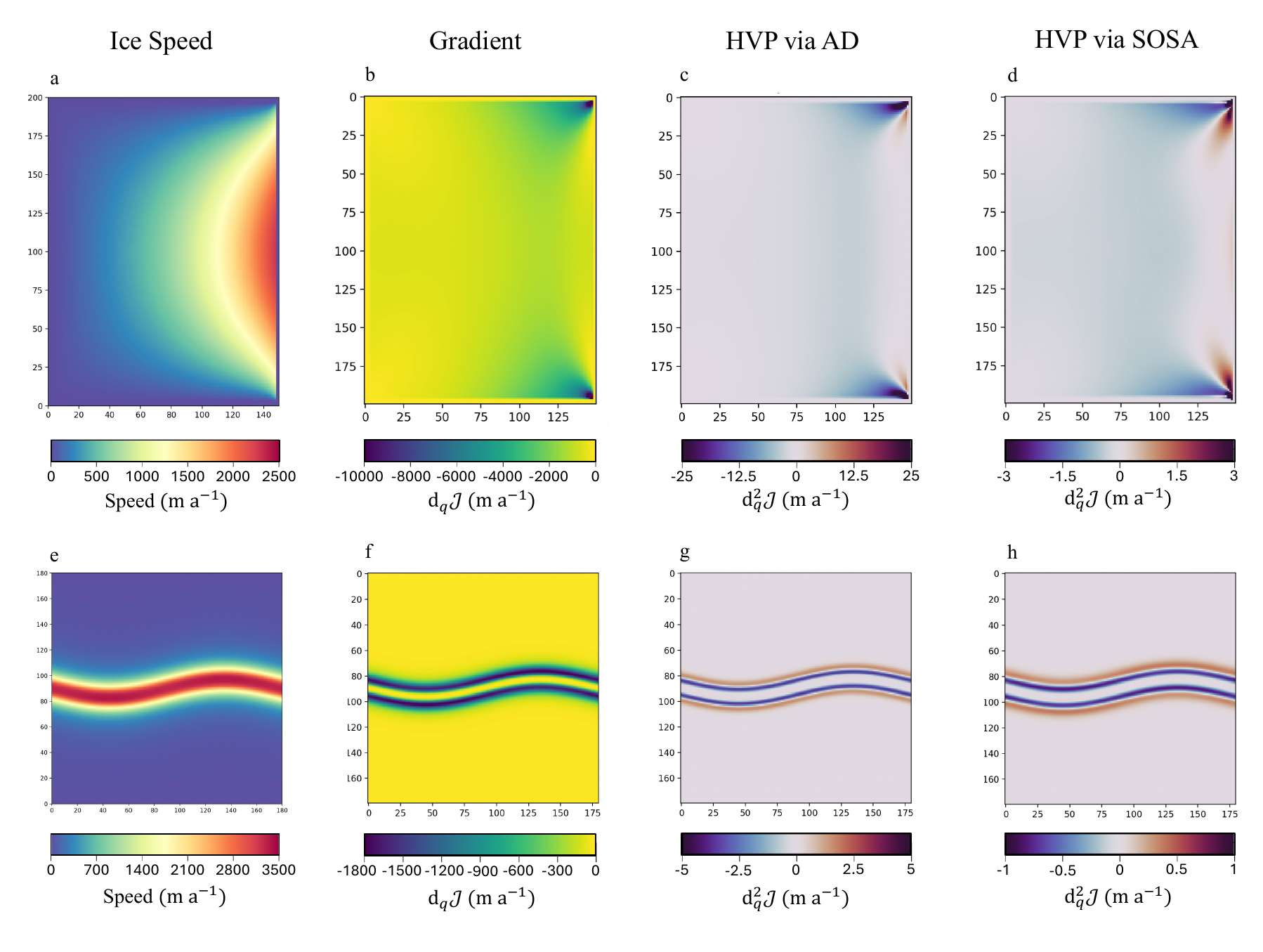}
\caption{Visual comparison of Hessian-vector products calculated using the AD and SOSA methods for the Ice Shelf (a-d) and Twisty Stream (e-h). The functional $\J$ measured the square-integrated ice speed over the domain. (a) and (e) show the speeds found by solving the SSA momentum balance equations. (b) and (f) show the the gradient $d_q\J$ computed via algorithmic differentiation (AD). (c) and (g) show the Hessian-vector product (HVP) $d^2_q\J[\dq]$ for a perturbation direction $\dq$ parallel to the gradient shown in (b) and (f) found using AD. (d) and (h) show the HVP calculated using the PDE-level second-order self-adjoint (SOSA) model.}\label{fig1}
\end{figure}

First, we show a visual comparison between example Hessian-vector products found using the SOSA and AD methods. We considered the functional $\J(\vecu) = \int_\Omega\sqrt{\vecu\cdot\vecu}~{\rm d}\Omega $. For both the Ice Shelf and Twisty Stream, we computed the ice speed, the gradient of $\J$ with respect to the parameter $q$ (defined in equation \ref{eq:ssa}), then the Hessian-vector product $d^2_q\J[\dq]$ where the perturbation $\dq= \epsilon d_q\J$ is in the gradient direction. Results are shown in Fig. \ref{fig1}.\\

In the ice shelf case, the gradient reflects the common observation that softening in the shear margins and in a central compressive band \cite{doake1998breakup} have largest impact on flow speed. The HVP in this direction, for both the AD and SOSA cases, display a similar structure. Both indicate a decreasing sensitivity to softening as it is applied upstream of the compressive arch, but an increasing sensitivity to it downstream. In the AD case, the magnitudes of the HVP are larger by up to an order of magnitude.\\

The Twisty Stream experiment shows that the ice speed, as expected, is most sensitive to softening in the shear margins of the ice stream. In both the AD and SOSA cases, the sensitivity to softening increases as it is applied at the outside edge of the current shear margins (corresponding to a widening of the ice stream and an increase in integrated flow speed), and a decreased sensitivity to softening at the inside edge (corresponding to a narrowing of the ice stream). Once more, this effect is amplified when the full adjoint is found using AD than in the SOSA case.\\

\subsection{The Hessian decomposition}

To compare the methods more concretely, we reduce our attention to the Twisty Stream example and consider the spectral decomposition of the Hessians computed via AD and SOSA. The domain is $180\times 180$-cells in size, so the Hessian has over 32 000 eigenvectors. We consider only the first 500 to get an idea of how the Hessians differ in their most consquential modes.\\

Given that we have in mind uses of the Hessian in the context of glaciological inverse problems, we look at a functional representing something like the misfit between modelled and observed ice speed:
\begin{equation}\label{cst_fct_ip}
    \J = \int_\Omega \sqrt{(\vecu-\vecu_{\rm sltn})\cdot(\vecu-\vecu_{\rm sltn})}~{\rm d}\Omega
\end{equation}

where $\vecu_{\rm sltn}$ is the ice velocity solution for unperturbed $q$. This has the feature of having zero gradient at $\vecu = \vecu_{\rm sltn}$.
\\

\begin{figure}[h!t]
\centering
\includegraphics[width=\linewidth,keepaspectratio]{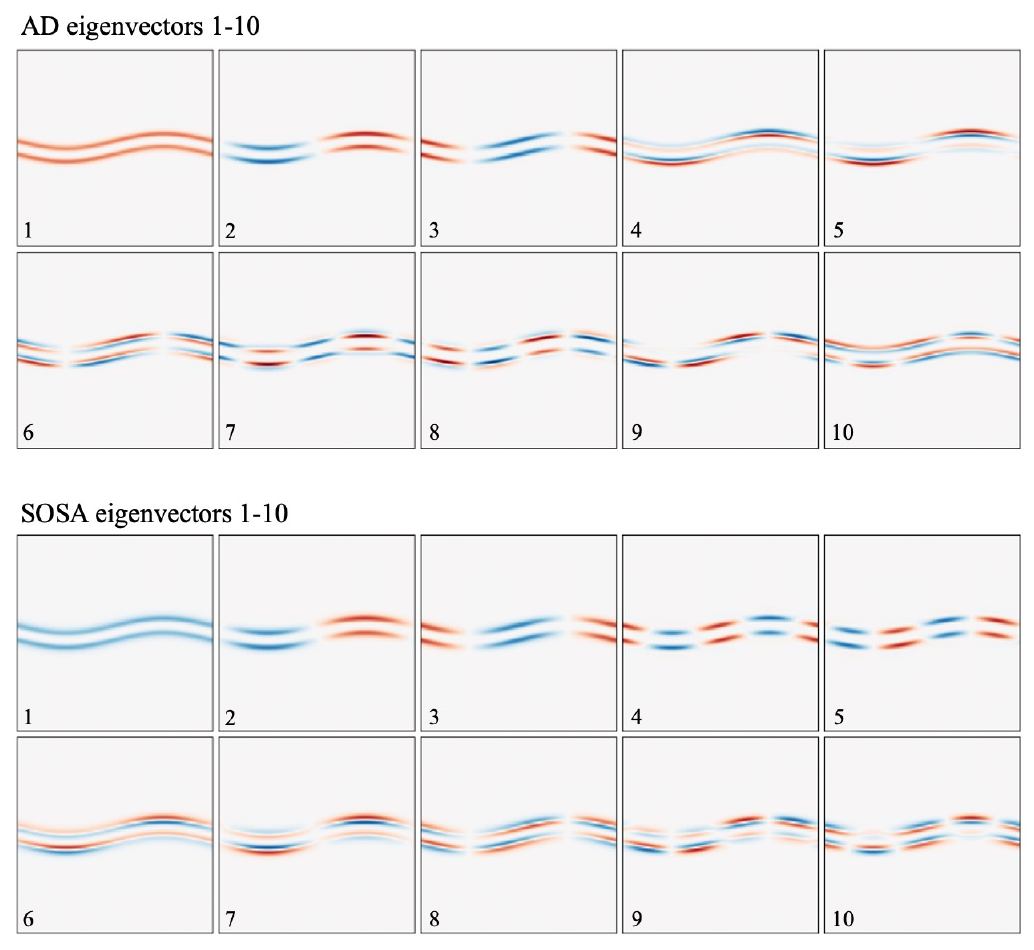}
\caption{The first 10 eigenvectors of the Hessian computed via the AD (upper) and SOSA (lower) methods. The functional in this case measures the difference in velocity from the solution at unperturbed $q$ (eq. \ref{cst_fct_ip}). Eigenvectors are ordered by the magnitudes of their eigenvalues from largest to smallest. The colourmap is centred on 0.}\label{evec_comp}
\end{figure}

Fig \ref{evec_comp} shows the first 10 eigenvectors produced using each method for this functional. Many similar patterns appear in the two sets, but their order is different, and often differ by the sign of the eigenvalue (e.g. the principal mode in each). For example, the pair of eigenvectors 4 and 5 of the AD set look similar to the 6-7 pair of SOSA eigenvectors. This pattern continues throughout the first 500 eigenvectors, with similar structures appearing in both sets. There is also a difference in the along-flow/across-flow periodicity in some early eigenvectors. For example, eigenvectors 4 and 5 of the SOSA-derived Hessian show half-wavelength in the across-flow direction and 2 wavelengths in the along-flow, and no such pattern appears in the first 500 eigenvectors of the AD-derived Hessian.\\

\begin{figure}[h!t]
\centering
\includegraphics[width=\linewidth,keepaspectratio]{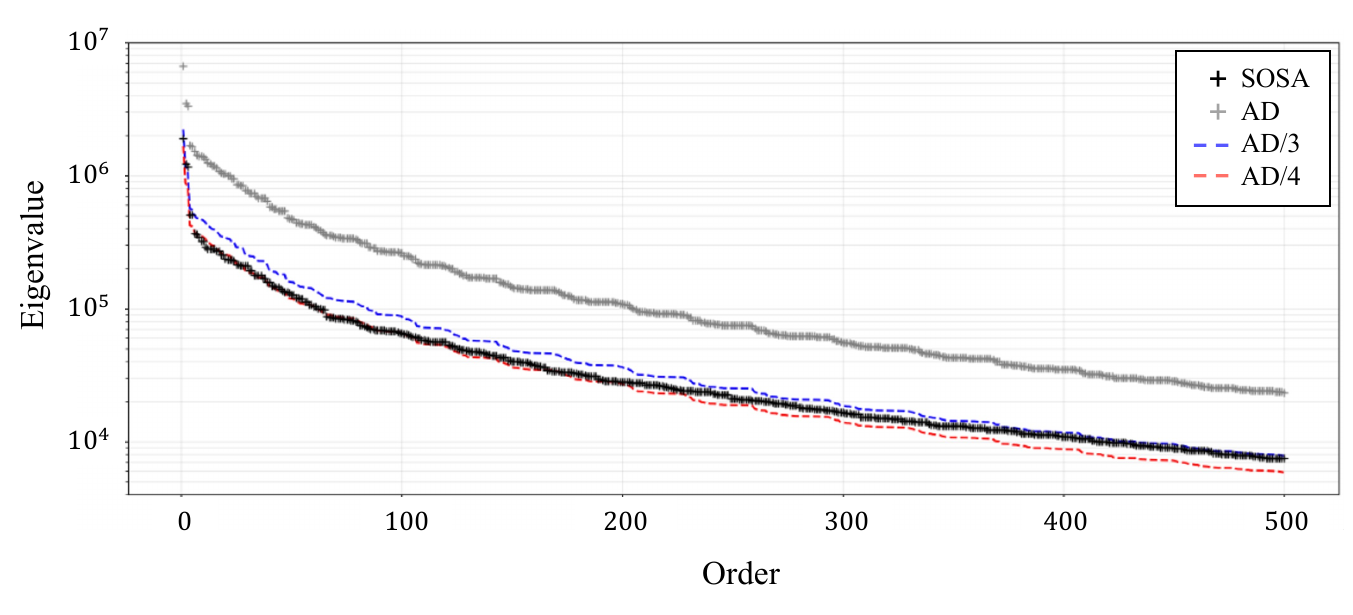}
\caption{The first 500 eigenvalues of the AD- and SOSA-derived Hessians, for the functional defined in \eqref{cst_fct_ip}, arranged in order of their magnitudes from largest to smallest.  Grey + icons show the AD eigenvalues, black + icons show the SOSA eigenvalues. The blue and red dashed lines show the AD eigenvalues divided by 3 and 4, respectively.}\label{evals}
\end{figure}

The eigenvalues for each Hessian reduce rapidly in magnitude for the first few modes, though go onto display sub-exponential decay (Fig. \ref{evals}). This behaviour has been seen in the highest order eigenvalues for similar problems before, and has been shown to give way to exponential decay as the order gets large \cite{bea_2023_paper}. The SOSA eigenvalues are smaller than the AD eigenvalues, due to the missing terms associated with derivatives of the viscosity. For the first few eigenvalues, the ratio is around a factor of 3 (Fig. \ref{evals} - blue line). This relates to the fact that the operator $\deu G$ is roughly a factor of $n=3$ (our choice for Glen's exponent) smaller when derivatives of the viscosity are included. This factor of $n$ propagates unchanged into the first term of equation \eqref{hvp_final} and into the eigenvalues for the eigenvectors in which this term dominates. After the first 5 modes, the difference drifts away from a constant factor, though remains largely between 3 and 4, though will drop below 3 eventually.\\

The similarity of the sets of eigenvectors under visual inspection means that greater analysis is needed to determine the utility of the SOSA approximation. For most applications in optimisation and uncertainty quantification, an approximation of the Hessian to some rank much smaller than the full dimension is necessary. For example, in an idealised optimisation procedure, one can imagine the search directions being described by the ordered eigenvectors of the Hessian of the cost function, and we don't need to proceed through all of them before being satisfied the problem is solved well-enough. So, we need to determine the similarity of the AD and SOSA Hessians as the rank of the approximation increases. Concretely, we ask at what point the subspaces spanned by the first k eigenvectors of each Hessian start to diverge.\\

\begin{figure}[h!t]
\centering
\includegraphics[width=\linewidth,keepaspectratio]{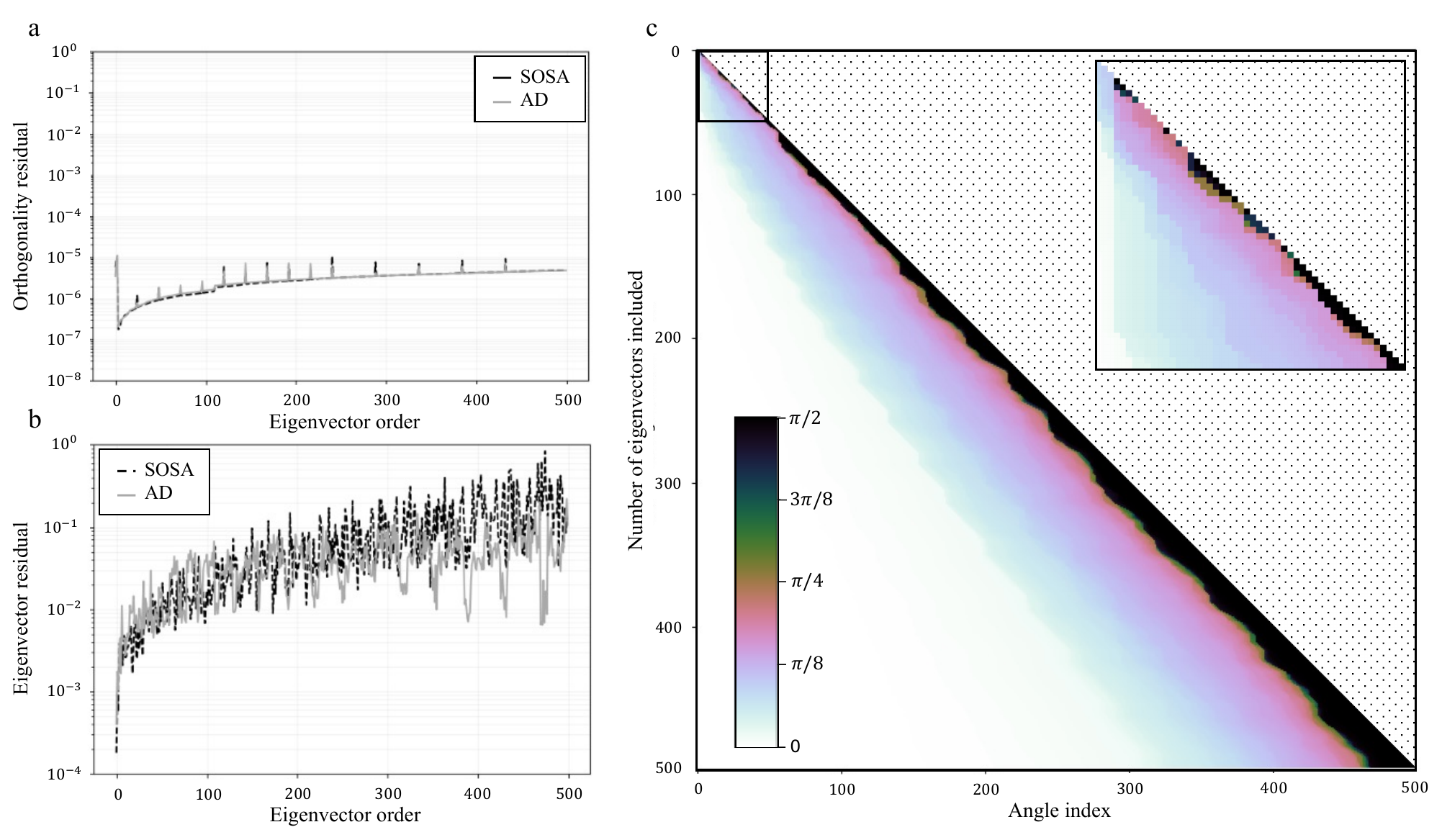}
\caption{The similarity of subspaces spanned by the first $k$ eigenvectors of the Hessian computed via the AD and SOSA methods. a) The residual \eqref{eigres} - indicating the quality of the eigevectors themselves. b) The residual \eqref{orthres} - indicating the self-consistency of each subspace. c) Ordered principal angles between the AD and SOSA subspaces for different numbers of eigenvectors. The black area indicates the angles close to $\pi/2$. Inset: principal angles between subspaces up to 50 eigenvectors.}\label{sms}
\end{figure}

To determine this, it is necessary to first consider whether the individual subspaces are themselves self-consistent, so that any divergence between the two reflects the true geometry rather than, for example, numerical error. One might worry about this, for example, because the methods include a number of linear algebra solves per HVP with poorly conditioned matrices. For each eigenvector/eigenvalue pair $(\vec{\hat{e}}_i, c_i)$, we computed an eigenvalue residual:
\begin{equation}\label{eigres}
{\rm eigres_i} = \frac{||d_q^2\J[\vec{\hat{e}}_i] - c_i\vec{\hat{e}}_i||_1}{||c_i\vec{\hat{e}}_i||_1}.
\end{equation}
Figure \ref{sms}b shows these remain small throughout the first 500 pairs computed for each method. Additionally, we tested the orthonormality of the subspaces by computing, for each subspace, the metric:
\begin{equation}\label{orthres}
{\rm orthres_k} = ||S^\top_k S_k - \mathcal{I}_k||_1,
\end{equation}
where $S_k = [\vec{\hat{e}}_1, \dots,\vec{\hat{e}}_k]$ and $\mathcal{I}_k$ is the $k\times k$ identity matrix. Figure \ref{sms}a shows that this too remains small up to $k=500$. The conclusion is that each method indeed produces a clean set of orthonormal eigenvectors.\\

An intuitive notion of the similarity of subspaces constructed from the leading modes of the AD and SOSA Hessians can be determined by considering the set of principal (or canonical) angles between them. Let $A$ be the matrix with columns defined by the first $k$ eigenvectors found using the SOSA method and $B$ be the equivalent matrix for the AD method. The singular value decomposition of the matrix $A^\top B$ defines the cosines of the principal angles between the subspaces defined by the columns of $A$ and $B$. When the angles are all close to zero, the projection of any vector in one subspace onto the other will be almost parallel to the first. When the largest of these angles is $\pi/2$, the subspaces are orthogonal, though the projection distance still might be small. As such, the proportion of angles significantly different from 0 gives us an indication of the similarity of the Hessians at different rank. Figure \ref{sms}c shows the ordered principal angles between subspaces constructed using the first $k=1\,{\rm to}\,500$ eigenvectors of the two Hessians.\\

Figure \ref{sms}c shows that the subspaces spanned by the eigenvectors of the SOSA and AD Hessians start out close to parallel but become persistently orthogonal after the 33rd mode. Up to the first 4 modes, the principal angles are all small ($<\pi/8$), however, after this, the proportion of larger angles grows until around the 50th mode (\ref{sms}c - inset). From here, the proportion of angles very close to 0 (white region) starts to dominate again. The number of angles near $\pi/2$ (black region) grows as modes are added, but it does not increase as a proportion of the full set of angles. Hence, though there is a persistent inaccessible region to the SOSA Hessian, increasing its rank appears to always improve the approximation. This general pattern implies that there is a set of ranks for which neglecting nonlinear terms makes the SOSA Hessian most different to the full Hessian: between $4$ and $\sim100$. This analysis is revealing about why gradient methods that assume self-adjointness appear to work well in practice, and has consequences for the use of the SOSA method in place of the full Hessian. We discuss these points below.

\section{Discussion}

\subsection{Optimisation}

As mentioned above, the context we have in mind for using the Hessian is the solution of inverse problems, for example for the purposes of model initialisation or diagnostic modelling, and quantifying uncertainty in the solutions. Commonly, solutions are found using gradient-based methods using only first order information. Despite not explicitly computing second-order derivatives, the geometry of the cost-landscape (e.g. \eqref{cst_fct_ip}) encoded in the Hessian still governs how those methods behave. We imagine, as an approximation, that reasonably close to the solution, steepest-descent-type methods (e.g. conjugate gradients) explore in turn the maximum-curvature directions corresponding to the eigenvectors of the Hessian. During an iterative optimisation algorithm using an approximate gradient, the solution will continue to improve only while the search direction remains non-orthogonal to the true gradient. The plot of principal angles between the AD- and SOSA-derived Hessians (Fig. \ref{sms}c) indicates that these directions are similar in the first few modes, then diverge slightly for the next few tens of modes, before divergence stops around mode 100. This gives us some insight into why gradient-based methods applying the self-adjoint approximation perform well in practice - there is no obvious point at which adding modes corresponding to inexact gradients stops improving the solution. Because they are expensive for real problems, optimisation algorithms are, in practice, run only for a few tens of iterations. However, Fig. \ref{sms} suggests improvement would continue on from this.\\

One potential use of the SOSA Hessian is in higher-order Newton-like methods, where we could use it to approximate the true Hessian at some low rank. The same considerations apply here as in the first-order case - in that the invariant subspaces spanned by the eigenvectors of each Hessian are similar enough to imagine this being effective. However, underestimation of curvature in each mode of the SOSA Hessian now becomes important. Depending on the rank to which the approximation is desired, a scaling factor could be applied to the eigenvalues (Fig. \ref{evals}). Given the difference in many of the leading-order curvature modes, it is unclear how well even a scaled SOSA Hessian would perform. Hence, where feasible, the full Hessian is preferable.\\

\subsection{Uncertainty Quantification}
Finally, in Bayesian parametric uncertainty quantification, the inverse Hessian is used to calculate the covariance of the posterior distribution for the parameter of interest \cite{thacker_roleOfHessian_1989, ISAAC2015348, Koziol2021}. The accuracy of the posterior covariance therefore places strict requirements on the fidelity of the Hessian. Given that the eigenvalues are different in the two cases, use of the SOSA Hessian would result in a systematic underestimation of uncertainty in each direction. To some degree, it seems that this could be compensated-for in the leading modes by applying a constant scaling to the inverse Hessian, though this would be quite inexact. Existing studies have shown the number of eigenvectors required for faithful construction of the posterior covariance numbers in the thousands, for problems not much larger than those considered here \cite{bea_2023_paper}. Even though the proportion of large angles between the SOSA and full Hessian might be small at this order, the differences at low order, and differences in eigenvalues, might well yield inaccurate results. However, there exist other methods of uncertainty quantification, such as perturbation ensembles, for which the broad alignment of the SOSA and full Hessian might suffice to produce similar results.\\

\section{Conclusion}

The increasing prevalence of algorithmic differentiation tools has led to the adoption of methods in ice sheet model initialisation that make use of second-order derivative information. The second-order adjoint model derived here, which neglects terms coming from the nonlinear dependency of ice viscosity on velocity, is an easily implemented alternative for the set of models not built on frameworks which support AD. The Hessian-vector products produced by this SOSA method, implemented in a JAX-based finite volume code, look, in the two cases we consider, to be visually similar, as do the first 500 eigenvectors of the Hessian in each case. Spectral analysis reveals that there is divergence between the subspaces spanned by eigenvectors of the full and approximate Hessians. This divergence is small enough to explain the success of first-order optimisation methods using the self-adjoint approximation. However, coupled with the difference in eigenvalues, it suggests use of the full Hessian, derived via AD, ought to be used where a high-fidelity approximation is sought at any rank higher than 3 or 4.

\vspace{0.5in}
\textbf{Acknowledgements}
TSS and SLC are supported by the Advanced Research and Invention Agency grant no. SCOP-PR01-P010.

\bibliographystyle{elsarticle-num}
\bibliography{refs}

\appendix
\section{First-order Adjoint Formulation}\label{appsec:foa}

We consider the sensitivities of a functional $\J(\vecu,q)\in \mathbb{R}$ with respect to changes in the control parameter $q$, for a model that solves a set of partial differential equations $G(\vecu,q)=0$. We take $q\in Q$ and $\vecu\in U$ where $Q$ and $U$ are both Hilbert spaces, so we use inner products throughout this derivation, and apply the Reisz representation theorem wherever convenient to associate gradients of functionals and Hessian-vector products with elements of those different spaces. \\

The ultimate quantity of interest in calculating the sensitivity of $\J$ changes in $\phi$, is the Gâteaux derivative of $\J$ in the direction $\delta q$:
\begin{equation}\label{gateaux_def}
    \delta\J(\vecu, q, \delta q) = \lim_{\varepsilon\rightarrow0}\left\{\frac{1}{\varepsilon}\left(\J(\vecu, q+\varepsilon\delta q) - \J(\vecu, q)\right)\right\}.
\end{equation}

We define the gradient $d_q\J$ as an operator, and as an element of $Q$ by equating all of the following things: 
\begin{equation}\label{eq:grad_direction_ip}
    \delta\J(\vecu, q, \delta q) = d_q\J[\dq] = \langle~d_{q}\J, \delta q~\rangle = \int_\Omega d_q\J(\vx)\, \dq(\vx)~d\Omega.
\end{equation}

The gradient $d_q\J$ can be decomposed via the chain rule, which operates in the same way in functional calculus as calculus:

\begin{equation}
    d_q\J = \deq\J + \deu\J\,d_q\vecu.
\end{equation}

The first-order adjoint method allows us to compute this without forming the Jacobian $d_q\vecu$ explicitly. There are a few ways of formulating it, but it is typical in the ice sheet modelling literature to imagine forming the Lagrangian:
\begin{equation}
    \mL = \J + \langle \vecl, G\rangle
\end{equation}
for some Lagrange multipliers $\vecl\in U$. These Lagrange multipliers are also known as the adjoint variables. The variation of this under perturbations of the control parameter $q$ is simply the variation of $\J$, given that $d_q G[\dq]=0$. (In other words, we still expect our solver to converge when $q$ is perturbed a small amount.) So, let us consider the action of the gradient of $\mL$ on a perturbation $\dq$:

\begin{equation}\label{lag_grad}
    d_q\J[\dq] \equiv  d_q\mL[\dq] = \big(\deq\J + (\deu\J\,d_q\vecu) + \deq\langle\vecl,G\rangle + (\deu\langle\vecl,G\rangle\,d_q\vecu)\big)[\dq].
\end{equation}

From the definition of the Gateaux derivative \eqref{gateaux_def}, it can be seen that:

\begin{equation}\label{obvs_identity_back}
    (\partial_{\vecu,q} \langle \vecl, G \rangle)[\dq] = \langle\vecl,\partial_{\vecu,q} G[\dq]\rangle,
\end{equation}
In which case, \eqref{lag_grad} gives us:

\begin{equation}
    d_q\J = \deq\J + (\deq G)^\ast[\vecl] + (\deu\J + (\deu G)^\ast[\vecl])d_q\vecu
\end{equation}
where we have used the definition of the adjoint of a linear map.\\

Hence, if we can solve the first-order adjoint equation:
\begin{equation}\label{generic_foa_eq}
    \deu\J + (\deu G)^\ast[\vecl] = 0
\end{equation}
for $\vecl$, then:
\begin{equation}\label{grad_given_foa}
    d_q\J = \deq\J + (\deq G)^\ast[\vecl].
\end{equation}

The question now becomes: what what are $(\deu G)^\ast[\vecl]$ and $(\deq G)^\ast[\vecl]$ for the SSA equations?

\subsection{The adjoint equation and gradient for the SSA}

Once, again the shallow-stream approximation equations are:
\begin{equation}
    \grad\cdot[\varphi(q)H(\vecu)] - C\vecu - \rho_i gh\grad s = 0.
\end{equation}
where $\varphi(q) = \phi_0e^q$ and $H(\vecu)$ is the resistive stress tensor, subject to boundary conditions:
\begin{equation}\label{eq:adj_bvp_2}
\begin{aligned}
    \hat{\vec{n}}\cdot\vecu = 0,\quad \hat{\vec{t}}\cdot\grad\vecu\cdot\hat{\vec{n}} = 0
\end{aligned}
\end{equation}
(where $\hat{\vec{n}}$ and $\hat{\vec{t}}$ are normal and tangent vectors to the boundary respectively).\\

The linear self-adjoint approximation, used widely in the literature, can be simply stated as:
\begin{equation}\label{linsa}
    (\deu G(\vecu))^\ast[\vecl] = G(\vecl).
\end{equation}

Given the structure of $G$, it is clear that making $\mub$ independent of $\vecu$ makes $G$ linear. In order to justify the self-adjointness, one must ensure that adjoint variables $\vecl$ also conform to the boundary conditions specified in \eqref{eq:adj_bvp_2}. The result of this is that all the boundary terms arising from integration by parts disappear when deriving \eqref{linsa}.\\

The $(\deq G)^\ast[\vecl]$ term is most easily computed by rewriting $(\deq G)^\ast[\vecl] = \deq\langle\vecl, G\rangle[\dq]$ and then direct application of the definition of the Gateaux derivative \eqref{gateaux_def}. Briefly:
\begin{align*}
    (\deq G)^\ast[\vecl] &= \deq\langle\vecl, G\rangle[\dq]\\
    &= - \lim_{\varepsilon\rightarrow0}\left\{\frac{1}{\varepsilon}\left(\int\nabla\vecl : ([\varphi(q+\varepsilon\dq) -\varphi(q)]H)~d\Omega)\right)\right\}\\
    &= -\nabla\vecl : (\varphi H).
\end{align*}
Where, the second line makes use of the following statement, which is just generalisation of Green's first identity: For an n-dimensional vector field $\vec{v}$ and a rank-2 tensor field $T$ defined over the domain $\Omega$:
\begin{equation*}
    \int_{\Omega}\vec{v}\cdot(\nabla\cdot T)~d\Omega = \int_{\partial\Omega}(T~\cdot\nhat)\vec{v}~d\Sigma - \int_{\Omega}T\colon(\nabla\vec{v})~d\Omega,
\end{equation*}
where $\nhat$ is the unit normal to the boundary $\partial\Omega$.\\

With this, we find:
\begin{equation}\label{eq:gradient}
    d_q\J \approx \deq\J - \grad\vecl:(\varphi(q)H(\vecu)).
\end{equation}

\section{Second-order Adjoint Formulation}\label{appsec:soa}

\subsection{General formulation}

We start by restating the adjoint equations we found in the derivation of the first-order sensitivity:
\begin{equation}\label{adjoint_eq_2}
    \partial_{\vecu}\J + \partial_{\vecu} G_{\vecl} = 0.
\end{equation}
These are still accompanied by the constraint that $\vecu$ solves the stress-balance equations:
\begin{equation}\label{constraint_2}
    G=0.
\end{equation}
The first step in calculating second-order derivatives is to solve the first order equations for the Lagrange multipliers $\vecl$. This allows us to write the functional gradient of the cost function as:
\begin{equation}
    d_q\J = \partial_{q}\J+\partial_{q} G_{\vecl}.
\end{equation}\\

Differentiating this and applying the chain rule gives us the Hessian:
\begin{equation}
    \mH = d^2_q\J = \de^2_q\J + \de^2_qG_{\vecl} + (\de_\vecu\de_q\J)d_q\vecu + (\de_\vecu\de_q\Gl)d_q\vecu + (\de_{\vecl}\de_q\Gl)d_q\vecl.
\end{equation}

It is instructive to consider the action of this Hessian on two perturbations $\dq_1$ and $\dq_2$:

\begin{equation}\label{Hessian_action}
\begin{aligned}
    \mH[\dq_1][\dq_2] = &(\de^2_q\J + \de^2_q\Gl)[\dq_1][\dq_2] ~+ \\
    &\langle~d_q\vecl[\dq_1],\de_q\de_{\vecl}\Gl[\dq_2]~\rangle ~+\\
    &\langle~d_q\vecu[\dq_1],(\de_q\de_\vecu\J + \de_q\de_\vecu\Gl)[\dq_2]~\rangle
\end{aligned}
\end{equation}

where we have used the fact that $(\partial_q\partial_\vecu\mathcal{F})^\ast = \partial_\vecu\partial_q\mathcal{F}$ for some functional $\mathcal{F}$ (as well as other similar identities).\\

As with the first-order adjoint system, derivatives of the constraints at our disposal are employed to reduce this to an equation that doesn't depend on Jacobians that would be prohibitive to compute. The new adjoint system will again constitute a system of equations whose solutions correspond to the Lagrange multipliers reqiured to do this.\\

To \eqref{Hessian_action}, we add two indetically zero terms: the inner product of two Lagrange multipliers $\vmu$ and $\vbeta$ with the gradients of the model and first-order adjoint equations applied to a perturbation $\dq_1$:

\begin{equation}
\begin{aligned}
    \mH[\dq_1][\dq_2] = &(\de^2_q\J + \de^2_q\Gl)[\dq_1][\dq_2] ~+ \\
    &\langle~d_q\vecl[\dq_1],\de_q\de_{\vecl}\Gl[\dq_2]~\rangle ~+\\
    &\langle~d_q\vecu[\dq_1],(\de_q\de_\vecu\J + \de_q\de_\vecu\Gl)[\dq_2]~\rangle ~+\\
    &\langle~\vbeta, d_q G[\dq_1]~\rangle ~+\\
    &\langle~\vmu, d_q(\deu\J + \deu\Gl)[\dq_1]~\rangle.
\end{aligned}
\end{equation}

After expanding and rearranging, we find:
\begin{equation}
\begin{aligned}
    \mH[\dq_1][\dq_2] = &(\de^2_q\J + \de^2_q\Gl)[\dq_1][\dq_2] ~+ \\
    &\langle~\vbeta, \deu G[\dq_1]~\rangle~+\\
    &\langle~\vmu, (\deq\deu\J + \deq\deu\Gl)[\dq_1]~\rangle~+\\
    &\langle~d_q\vecl[\dq_1],\de_q\de_{\vecl}\Gl[\dq_2] + \deu\del\Gl[\vmu]~\rangle ~+\\
    &\langle~d_q\vecu[\dq_1],(\de_q\de_\vecu\J + \de_q\de_\vecu\Gl)[\dq_2] + (\deu G)^\ast[\vbeta] + (\deu^2\Gl + \deu^2\J)[\vmu]~\rangle.
\end{aligned}
\end{equation}

This shows that we can formulate a Hessian-vector product of the form:

\begin{align}\label{eq:svp_expansion}
    \mH[\dq] = (\de^2_q\J + \de^2_qG_{\vecl})[\dq]+(\de_\vecu\de_q\J + \de_\vecu\de_q\Gl)[\vmu] + (\de_q G)^\ast[\vbeta]
\end{align}

if we choose the $\vmu$ and $\vbeta$ that solve the equations:

\begin{align}
    &\de_q\de_{\vecl}\Gl[\dq_2] + \deu\del\Gl[\vmu] = 0\label{mu_eq}\\
    &(\de_q\de_\vecu\J + \de_q\de_\vecu\Gl)[\dq_2] + (\deu G)^\ast[\vbeta] + (\deu^2\Gl + \deu^2\J)[\vmu] = 0.
\end{align}

These are the second-order adjoint equations.

\subsection{Application to the SSA}

To apply this to the shallow-stream approximation to the Cauchy momentum equations for ice sheet flow, there are 12 terms in the above equations that we need to evaluate. Choosing a cost functional that is the sum of separate functionals of $\vecu$ and $q$ means that terms involving $\de_\vecu\de_q\J$ can be set to zero. As with the calculation of first-order derivatives, we assume a viscosity that is independent of $\vecu$, so that $\mub$ doesn't depend on $\vecu$ or $q$. This means that $G$ is linear in $\vecu$ and we can set $\partial_\vecu^2 G_{\vecl} =0$ as well. This leaves five expressions in the SOA system (now the ``SOSA'' system, after the afforementioned approximation) that need to be evaluated, and four in the equation for the Hessian-vector product. These are:

\begin{align}
    &\de_q\de_{\vecl}\Gl[\delta q]\label{soat1}\\
    &\de_{\vecu}\de_{\vecl}\Gl[\vmu]\label{soat2}\\
    &\de_q\de_{\vecu}\Gl[\delta q]\label{soat3}\\
    &\de_{\vecu}^2\J[\vmu]\label{soat4}\\
    &(\de_{\vecu}G)^\ast[\vbeta]\label{soat5}\\
    &\de_q^2\J[\delta q]\label{hvp6}\\
    &\de_q^2\Gl[\delta q]\label{hvp7}\\
    &\de_q\de_\vecu\Gl[\vmu]\label{hvp8}\\
    &(\de_q G)^\ast[\vbeta]\label{hvp9}.
\end{align}

We will deal with these in order, starting with the first 5 which define the SOSA system:

\begin{itemize}
    \item {
    Eq. \eqref{soat1}.
    \begin{align}
        \deq\del\Gl[\dq] &= \deq G[\dq]\\
        &= \lim_{\varepsilon\rightarrow0}\left\{\frac{1}{\varepsilon}\left(G(q+\varepsilon\delta q, \vecu) - G(q, \vecu)\right)\right\}\\
        &= \nabla\cdot (\varphi~\dq~ H(\vecu))
    \end{align}
    where we have used that everything is smooth and linear.}

\item{
Eq. \eqref{soat2}:

\begin{align}
    \deu\del\Gl[\vmu] &= \deu G[\vmu] = G(\vmu)\\
    &= \nabla\cdot(\varphi~ H(\vmu))-C\vmu.
\end{align}
This makes the crucial SOSA assumption that $G$ is linear.}

\item{
Eq. \eqref{soat3}:

Consider the action of $\deu \Gl$ on $\delta\vecu$. As $\Gl$ is linear, this is just $\Gl(\delta\vecu)=\int\vecl\cdot\left[\nabla\cdot(\varphi~\delta q~ H(\delta\vecu)) - C\delta\vecu\right]~d\Omega$. Consider:
\begin{align}
    \langle~ \deq \Gl(\delta\vecu), \delta q ~\rangle &= \int\vecl\cdot\nabla\cdot(\varphi~\delta q~ H(\delta\vecu)) d\Omega\\
    &= \int\delta\vecu\cdot\nabla\cdot(\varphi~\delta q~ H(\vecl)) d\Omega \label{final_but}
\end{align}
where the second line is once more due to the self-adjointness of the viscous part of $G$. Eq. \eqref{final_but} shows the action of the mixed derivative $\deq\deu\Gl$ on $(\delta q, \delta\vecu)$. Hence the map $\deq\deu\Gl[\delta q]$ can be seen to take the form:
\begin{equation}
    \deq\deu\Gl[\delta q] = \nabla\cdot(\varphi\delta q H(\vecl)).
\end{equation}
}

\item{
Eq \eqref{soat4}: We assume that, by construction, $\deu^2\J[\vmu]$ is not too difficult to write down.
}

\item{

Eq. \eqref{soat5}:
\begin{align}
    (\deu G)^\ast[\vbeta] &= \nabla\cdot(\varphi H(\vbeta)) - C\vbeta.
\end{align}
This is the self-adjointness again.
}
\end{itemize}

With these in the bag, the second-order adjoint equations for the shallow-stream approximation to the Stokes equations are:

\begin{align}\label{soa_system}
    &\grad\cdot(\varphi\dq H(\vecu)) + \nabla\cdot(\varphi H(\vmu)) - C\vmu = 0\\
    &\grad\cdot(\varphi\dq H(\vecl)) + \nabla\cdot(\varphi H(\vbeta)) - C\vbeta + \deu^2\J[\vmu]= 0.
\end{align}

For the case in which $\J = \int(\vecu-\vecu_{\rm obs})^2~d\Omega$, the term $\deu^2\J[\vmu]$ is equal to $2\vmu$.\\

Once again, the derivation of these equations has required the we set a litany of boundary terms to zero. Hence we require all variables obey the same boundary conditions as $\vecu$\\

Once these equations are solved, the Hessian-vector product can be evaluated according to eq. \eqref{eq:svp_expansion}. Let us derive \eqref{hvp6}-\eqref{hvp9}.

\begin{itemize}
    \item{Eq \eqref{hvp6}. $\deq^2\J[\dq]$ will depend on the form of the regularisation used in the cost function.}

\item{
Eq \eqref{hvp7}:
\begin{align}
    \deq^2\Gl[\dq] &= \deq\langle\deq\Gl,\delta q\rangle\\
    &=\deq\int\vecl\cdot\nabla\cdot(\varphi~\delta q~ H(\vecu))~d\Omega\\
    &= -\de_q\int(\nabla\vecl:(\varphi~\delta q~ H(\vecu)))\\
    &= -\nabla\vecl:(\varphi~\delta q~ H(\vecu))
\end{align}}

\item{
The same kind of reasoning gives:
Eq. \eqref{hvp8}:
\begin{align}
    \deu\deq\Gl[\vmu] = -\nabla\vecl:\varphi H(\vmu)
\end{align}}

\item{
Eq. \eqref{hvp9}:
\begin{align}
    (\deq G)^\ast[\vbeta] = -\nabla\vbeta:\varphi H(\vecu).
\end{align}}

\end{itemize}

Hence, the Hessian-vector product can be calculated according to:

\begin{align}\label{eq:hvp_final}
    \mH[\dq] = \de^2_q\J[\delta q] - \nabla\vecl:(\varphi~\delta q~ H(\vecu)) -\nabla\vecl:\varphi H(\vmu) -\nabla\vbeta:\varphi H(\vecu).
\end{align}

Taking a linear rheology for the calculation of first and second order sensitivities has led to a system of equations for the adjoint variables $\vecl$, $\vmu$, $\vbeta$ that all essentially have the same form, with different right-hand sides. This means little development is required to perform the calculations necessary for approximating Hessian-vector products once a forward SSA code has been set up. This makes the SOSA model derived here an appealing drop-in to many current ice sheet models, requiring only the use of existing stencils.

\end{document}